\documentclass{article}
\usepackage{amsmath,amsrefs,amsfonts}
  \usepackage{paralist}
  \usepackage{graphics} %% add this and next lines if pictures should be in esp format
  \usepackage{epsfig} %For pictures: screened artwork should be set up with an 85 or 100 line screen
 \usepackage[colorlinks=true]{hyperref}
\hypersetup{urlcolor=blue, citecolor=red}

\newtheorem{theorem}{Theorem}[section]

\newtheorem{remark}{Remark}
\usepackage{graphicx,booktabs}
\usepackage{mathrsfs}
\usepackage[linewidth=1pt,style=1]{mdframed}
\usepackage{caption,subfig}

\newcommand{\beq}{\begin{equation}}
\newcommand{\eeq}{\end{equation}}

\def\H01{H_0^1}

\def\yd{y^\delta}

\def\M{{\mathcal{M}}}
\def\<{\langle}
\def\>{\rangle}
\newcommand{\cz}{\ensuremath{C_0}-semigroup }
\def\d{\partial}

\newcommand{\argmin}[1]{\ensuremath{\mbox{arg}\,\, \underset{#1}\min}}

\newcommand{\rst}[1]{\ensuremath{{\mathbin\mid}%\upharpoonright}%
_{\raise-1.5ex\hbox{$#1$}}}}

\begin{document}
\title{An adjoint control method for initial condition identification of the Abstract Cauchy problem}
\author{}
\maketitle
%% Enter the first author's name and address:
\centerline{\scshape Cary Humber }
\medskip
{\footnotesize
 %% please put the address of the first author
 \centerline{Naval Surface Warfare Center}
   \centerline{Panama City, FL, 32407, USA}
} %% Do not forget to end the {\footnotesize by the sign }

\medskip

\centerline{\scshape Kazufumi Ito }
\medskip
{\footnotesize
 %% please put the address of the second author
 \centerline{Department of Mathematics, North Carolina State University}
   \centerline{Raleigh, NC 27695, USA}
} %

\bigskip

%% The name of the associate editor will be entered by an editorial staff
 \centerline{(Communicated by the associate editor name)}
%\author{Cary Humber$^{\rm a\ast}$\thanks{$^\ast$Corresponding author. Email: cary.humber@navy.mil
%\vspace{6pt}} and Kazufumi Ito$^{\rm b}$\\\vspace{6pt}  $^{\rm a}${\em{Naval Surface Warfare Center, Panama City, FL}}; 
%$^{\rm b}${\em{Center for Research in Scientific Computation, North Carolina State University, Raleigh, NC}}\\\vspace{6pt}
%\received{v3.4 released May} }

\begin{abstract}
\vspace{2mm}

This paper develops and analyzes a generic method for reconstructing solutions to the abstract Cauchy problem in a general Hilbert space, from noisy measured data.
The method is based on the relationship between a partial differential equation and its adjoint equation with control.  We demonsrate the capability of the method through 
analysis and numerical experiments.
\end{abstract}
\maketitle
\section{Problem Description}

%We are interested in problems of the form
%Inverse problems constitute an important category of problems in the sciences.  
%In this paper, we consider the abstract Cauchy problem
Let $X$ be a Hilbert space endowed with the inner product $\<\cdot,\cdot\>_X=\|\cdot\|_X^2$ and let $A:D(A)\subset X\to X$ be the infinitesimal generator
of a strongly continuous semigroup $S_t:\mathbb{R}\to\mathcal{L}(X)$.  
We are concerned with the following problem.  Given a $C^1$ function $f:[0,t_f]\to X$ determine the initial condition, $x_0$,
of the Cauchy problem
\begin{equation}
\frac{dx}{dt}(t)  =  Ax(t) +f(t)
\label{cauchy}
\end{equation}
satisfying the measurement condition
%where we obtain measurements
\begin{equation}
y(t) =  Cx(t) .
\label{meas}
\end{equation}
%for $x$ in a Hilbert space $X$ (e.g., $X=L^2(\Omega)$).  We assume the operator $A$ is linear and generates a strongly continuous
%semigroup $S_t:\mathbb{R}\to\mathcal{L}(X)$, 
The bounded operator $C:X\to Y$ retains information about the solution, which may only be a portion of the solution, such as
boundary values.
It is desired to determine the initial condition $x_0$ given the incomplete (and possibly noisy) measurements $y(t),\,\, 0\le t\le t_f.$
The method developed in this paper is capable of forecasting future states as well, however, we focus on the inverse problem of determining
$x(0)$ due to its practicality and the necessity of dealing with this more difficult problem.
We are especially interested in the case of partial measurements 
(i.e., the measurements are sparsely distributed over the domain $\Omega$).
In the following section, we develop methods for determining the initial state $x(0)$ and we demonstrate how the same methods are applicable to 
forecasting the state $x(t_f)$ with only minor adaptations.
We assume the Hilbert space $X$ is separable, so that there exists a complete orthonormal sequence $\{\varphi_k\}_{k=0}^\infty$ in $X$.  
The approximation of $x(0)=x_0$ is given by the truncated (generalized) Fourier series
\begin{equation*}
x_0^m=\sum\limits_{k=0}^m\alpha_k\varphi_k,
\end{equation*}
where the coefficients satisfy $\alpha_k=\<x_0,\varphi_k\>$.  Thus, within this framework, the problem reduces to estimating
the generalized Fourier coefficients of $x_0$.
The problem of identifying the initial condition of the abstract Cauchy problem has been widely studied.
Methods concerning this problem have been covered by Auroux and Blum \cite{bfn}, Ito \textit{et al} \cite{timereverse}, and references therein.
The monograph by Isakov \cite{isakov2006inverse} covers inverse problems for PDE in detail.  
%Other methods that have been studied are regularized linear least squares methods, ...
As with many inverse problems, there is extreme difficulty in recovering a function from partial and noisy measurements,
thus suitable regularization is necessary.  The main focus of this paper is on the reconstruction method coupled with
the multi-parameter Tikhonov regularization \cite{newchoice,regparam}.

\section{An Adjoint Method for approximating the Fourier expansion of the initial condition}
\label{dualcontrol}

In this section, we develop and analyze a new approach for estimating the initial condition of the abstract Cauchy problem
\eqref{cauchy} from time-series data.   The method developed here involves an indirect
computation of the generalized Fourier coefficients, based on the adjoint equation of the Cauchy problem.
This method has a direct link with optimal control theory.
Given noisy data, the accuracy and stability of the method will be demonstrated in this paper.
The general framework of our method allows any PDE formulated under the linear semigroup theory to fit into this framework.  \
Not only that, but it will be shown that 
the method can also be applied to the less ill-posed problem of forecasting future states of the system.  Thus, our method may be especially 
beneficial for applications such as weather forecasting or financial futures, where it may be necessary to go both backward and forward.

Consider the adjoint equation of \eqref{cauchy}, given by
\begin{equation}
-\frac{d p}{d t}(t) = A^*p(t)+C^*u(t) \quad %, \quad p(T)=0. 
\label{d2}
\end{equation}
where $C^*\in\mathcal{L}(Y,X)$ corresponds to the adjoint of the observation operator $C$, and, likewise, 
$A^*$ is the adjoint of the infinitesimal generator $A$.
Here, $u\in L^2(0,t_f;Y)$ denotes a control or input to the system.  
It will be demonstrated that a suitable control 
can be determined for which the generalized Fourier coefficients can be approximated by a combination of the control, $u$, and the data, $y$.

Recall the state equation is given by
\begin{equation}
\frac{dx}{dt}(t)=Ax(t) + f(t)
\label{d1}
\end{equation}
and the measurements, satisfying
\begin{equation}
y(t)=Cx(t),\quad 0\le t\le t_f
\end{equation}
are given, for a known source $f$.
Multiplying (\ref{d1}) by \(p\), (\ref{d2}) by \(x\), subtracting and integrating over \((0,t_f)\) yields
\begin{eqnarray}
 \int\limits_0^{t_f}\frac{d}{dt}\< x(t),p(t)\>\, dt = \int\limits_0^{t_f}\left(\< Ax,p\> - \< A^*p,x\> - \< C^*u,x\> + \< f(t),p(t)\>\right)\, dt
\label{}
\end{eqnarray}
which implies that
\begin{equation}
\< x(t_f),p(t_f)\>_X - \< x(0),p(0)\>_X =  \int\limits_0^{t_f} \< f(t),p(t)\>_X-\< u(t),Cx(t)\>_Y\, dt
\label{}
\end{equation}
yielding the relation
\begin{equation}
\<x(t_f),p(t_f)\>_X-\< x(0),p(0)\>_X = \int\limits_0^{t_f}\< u(t),\xi(t)-y(t)\>_Y  dt,
\label{relation}
\end{equation}
where
\begin{equation*}
\xi(t)=\int_0^tC S_{t-s}f(s)\, ds.
\end{equation*}
The relationship \eqref{relation} forms the foundation for approximating $\<x_0,\varphi_k\>$.
% (or $\<x_{t_f},\varphi_k\>$).

We recall that the unique mild solutions of the abstract Cauchy problem and its dual, with conditions $x(0)=x_0,p(t_f)=p_{t_f},$
are respectively given by
\begin{equation}
x(t) = S_tx_0 + \int_0^t S_{t-s}f(s)\, ds ; \quad 
p(t) = S^*_{t_f-t}p_{t_f}+\int\limits_t^{t_f}S^*_{s-t}C^*u(s)\, ds
\end{equation}  
for each $t\in[0,t_f].$
For reconstructing the initial state $x_0$,
we assume the controllability of the adjoint system \eqref{d2}, which is equivalent to the observability of
\eqref{cauchy}-\eqref{meas} (i.e., the pair $(A,C)$ is observable).   The pair $(A,C)$ is observable if for all $x\in X$
\beq
\int_0^{t_f}\| CS_t x\|_Y^2\, dt\geq \gamma\| x\|^2
\label{null}
\eeq
for some $\gamma > 0$.
%Note that the functions \(x,u,p\) are allowed to be vector valued.
%Now, suppose \(p(0)\) is known and define the linear operator \(\mathscr{L}\) by
Note that by the observability assumption \eqref{null}, the equation
\begin{equation*}
y=\M x
\end{equation*} admits a unique solution for $y \in R(\M)$, where $\M$ is the operator defined by 
\begin{equation*}
\M:= CS_t\quad 0\le t\le t_f.
\end{equation*}
  Furthermore, this unique solution depends
continuously on $y$.
The details of infinite dimensional control theory are covered in \cite{zwart}.
Having assumed the controllability of \eqref{d2}, we define the operator $\mathscr{L}:L^2(0,t_f;Y)\to X$ by
\begin{equation}
\mathscr{L} u := \int\limits_0^{t_f} S^*_{s}C^*u(s)\, ds
\label{control}
\end{equation}
for $u\in L^2(0,t_f;Y).$
By construction, the adjoint equation \eqref{d2} evolves backwards in time.
If $p(t_f)=0$, then the adjoint satisfies  
%By setting $p_{t_f}=0$, we obtain
$$
p(0)=\int_0^{t_f} S^*_sC^*u(s)\, ds.
$$

%and seek a solution to
%which means that $p(0)\in R(\mathscr{L}).$
By the controllability/observability assumption, we know a unique solution to 
\begin{equation}
\mathscr{L} u=p(0),
\label{control1}
\end{equation}
%\eqref{control1}, $u,$ 
exists for $p(0)\in R(\mathscr{L}).$
However, in practice, the exact controllability of \eqref{d2} is, in general, not true, so we assume the condition \eqref{control} holds approximately, i.e., there exists $u_\varepsilon$ such that
\beq
\| \mathscr{L} u_\varepsilon-p(0)\|_X\le\varepsilon
\label{controlapprox}
\eeq
for any $\varepsilon>0$.
Whenever this relationship does not hold (or holds only approximately), we must suitably regularize the problem, so that a reasonable $u$ can be obtained.
We are interested in solving for $u$, since by relationship \eqref{relation}, we can obtain the $k^{\rm th}$ Fourier coefficient of $x_0$ as 
$$
\<x_0,\varphi_k\>_X=\int_0^{t_f}\<u(t),y(t)-\xi(t)\>_Y\, dt,
$$
whenever $p_{t_f}=0.$

We proceed by defining a collection of adjoint functions $p_k(0)=\varphi_k$, such that $\{\varphi_k\}_{k=0}^m$ 
forms an orthonormal basis for a finite-dimensional subspace $X_m\subset X$.
Then $\{\<x(0),\varphi_k\>\}_{k=0}^m$ are the generalized Fourier coefficients for $x(0)$.  
By the controllability assumption \eqref{null} and by utilizing relation \eqref{relation}, we can determine the Fourier coefficients of $x(0)$ by solving the operator equations
\begin{equation}
\mathscr{L} u_k =\varphi_k,\quad 0\le k \le m.
\label{condition}
\end{equation}
If \eqref{control} or \eqref{controlapprox} holds, we will construct stable approximations $u_k$ using a suitable 
regularization method.    An example of such a regularization method for determining one-dimensional $u_k$ is to solve the minimization problem
\begin{equation*}
\min\limits_{u\in L^2(0,t_f;Y)} \|\mathscr{L} u-\varphi\|_X^2+\eta_1\int_0^{t_f}|u(t)|\, dt + \frac{\eta_2}{2}\int_0^{t_f}|u'(t)|^2\, dt,
\end{equation*}
where the first term corresponds to the sparsity of the approximate solution $u_k(t), t\in [0,t_f],$ while the second term corresponds to the smoothness of $u_k$.
We note that the smoothness of $u_k$ may affect noise dampening (see Remark \ref{noiseremark}).  
Such regularization methods are described in detail in the papers \cite{newchoice}, along with criteria for selecting the regularization parameters $\eta_1,\eta_2$.

Our approach is based on the fact that for each basis function $\varphi_k$ there exists a control $u_k\in L^2(0,t_f;Y)$
such that $\mathscr{L} u_k=\varphi_k$ (or $\| \mathscr{L} u_\varepsilon-p(0)\|_X\le\varepsilon$).  The controls $u_k(t)$ are determined in such a way that each 
adjoint $p_k$ is driven from zero at time $t_f$ to $p_k(0)=\varphi_k$, for a suitably chosen $\varphi_k$. With each $u_k$ determined, we 
construct the approximation for $x_0$ by
\begin{equation*}
x_0^m=\sum\limits_{k=0}^m\alpha_k\varphi_k
\end{equation*}
where the generalized Fourier coefficients are approximated by
\beq
\<x_0,\varphi_k\>_X\approx\int_0^{t_f} \<u_k(t),y(t)-\xi(t)\>_Y\, dt=\alpha_k
\label{approxfourier}
\eeq
using the relation \eqref{relation} and equation \eqref{condition}.

Further analysis of the method is detailed below, including the error analysis in Theorem \ref{errortheorem}.
The following summarizes the method for estimating $x_0$.
\begin{mdframed}[leftmargin=10pt,rightmargin=10pt]
\textbf{Dual Method for reconstruction of $x_0$:}\label{dualalg}
\begin{enumerate}
\item Pick an orthonormal basis, $\{\varphi_k\}_{k=0}^m$ for $X_m\subset X$
\item For each $k$ solve $\mathscr{L}u_k=\varphi_k$ to find  $u_k\in L^2(0,t_f;Y)$
\item Form the estimate for $x_0,$ \begin{equation*} x_0^m=\sum\limits_{k=0}^m\alpha_k \varphi_k\end{equation*}
where 
\begin{equation*} \alpha_k=\int_0^{t_f}\< u_k(t),y(t)-\xi(t)\>_Y\, dt
\end{equation*}
with
\begin{equation*}
\xi(t)=\int_0^{t}CS_{t-s}f(s)\, ds.
\end{equation*}
\end{enumerate}
\end{mdframed}

The well-posedness of the method follows from the controllability assumption \eqref{null}. 
%i.e. there exists $\gamma >0$ such that
\beq
\int_0^{t_f} \| CS_tx\|^2\, dt \ge \gamma \| x\|_X^2
\label{exactcontrol}
\eeq
for all $x\in X.$

\subsubsection*{\underline{Using the method for forecasting a future state}}
$ $\newline

\noindent
Now, we briefly introduce how the method is utilized for the purpose of forecasting a future state $x(t_f)$.
For this purpose, we assume the adjoint \eqref{d2} is null-controllable, i.e. there exists $u\in L^2(0,t_f;Y)$ such that $p(0)=0$ and
\beq
\mathscr{L} u=- S_{t_f}p(t_f).
\label{finalcondition}
\eeq
Recall that $p$ evolves backwards in time (with respect to the evolution of $x$).
In general, the exact null-controllability may not hold, however we assume the condition \eqref{finalcondition}
holds approximately, i.e., there exists $u_\varepsilon$ such that
\begin{equation*}
\|\mathscr{L} u_\varepsilon+S_{t_f}p(t_f)\|_X\le\varepsilon,
\end{equation*}
for any $\varepsilon>0.$
With $u_k$ determined, the generalized Fourier coefficients are approximated by
\begin{equation*}
\<x(t_f),\varphi_k\>_X=\int_0^{t_f}\<u_k(t),\xi(t)-y(t)\>_Y\, dt
\end{equation*}
where $u_k$ is the approximate solution to
\begin{equation*}
\mathscr{L} u_k=- S_{t_f}\varphi_k.
\end{equation*}
For the final state case, the method is well-posed under the assumption of 
null-controllability of the adjoint control system, i.e.
\begin{equation*}
S^*_{t_f}X\subseteq R(\mathscr{L}).
\end{equation*}
The method is summarized as follows:
\begin{mdframed}[leftmargin=10pt,rightmargin=10pt]
\textbf{Dual Method for reconstruction of $x_{t_f}$:}\label{dualalgfinal}
\begin{enumerate}
\item Pick an orthonormal basis, $\{\varphi_k\}_{k=0}^m$ for $X_m\subset X$
\item For each $k$ solve $\mathscr{L}u_k=-S_{t_f}\varphi_k$ to find $u_k\in L^2(0,t_f;Y)$
\item Form the estimate for $x_{t_f},$ \begin{equation*} x_{t_f}^m=\sum\limits_{k=0}^m\alpha_k \varphi_k\end{equation*}
where 
\begin{equation*} \alpha_k=\int_0^{t_f}\< u_k(t),\xi(t)-y(t)\>_Y\, dt
\end{equation*}
with
\begin{equation*}
\xi(t)=\int_0^{t}CS_{t-s}f(s)\, ds.
\end{equation*}
\end{enumerate}
\end{mdframed}

The novelty of this method is, in part, due to the fact that it is not necessary to compute
the time history of the adjoint, $p$.  However, the method utilizes the information available from the adjoint in order
to accurately reconstruct $x_0$.  
By utilizing the $L^1$ norm, we are able to construct sparsely distributed controls,
which can aid computational efficiency.
Furthermore, the method is quite robust to noise, as the actual inverse problem
does not involve the noisy data.  
\begin{remark}
We also note that there is a stochastic interpretation of this method.
Assume $x,p$ are random variables satisfying the linear stochastic differential equations
\begin{equation}
dx =(Ax(t)+f(t))dt+\sigma dB_t;\quad
-dp =(A^*p(t) +C^*u(t))dt 
\end{equation}
where $B_t$ is the Brownian motion, and $\sigma$ is the standard deviation (diffusion coefficient).
Then, by the relation \eqref{relation} we have
\begin{equation*}
\<x_0,\varphi_k\>_X=\int_0^{t_f}\<u_k(t),y(t)\>_Ydt-\int_0^{t_f}\<f(t),p_k(t)\>_X dt +\sigma\int_0^{t_f} p_k(t) dB_t
\end{equation*}
which implies that 
\begin{equation*}
E[|\<x_0,\varphi_k\>_X-\int_0^{t_f}\<u_k(t),y(t)-\xi(t)\>_Ydt|^2]=E[\sigma^2|\int_0^{t_f}p_k(t)\, dt|^2].
\end{equation*}
Thus, the mean square error in approximating the Fourier coefficients is proportional to the standard deviation, $\sigma$, of the Brownian motion, 
regardless of that fact that $p(t_f)=0$ (in the case of estimating $x_0$).
Determining the control, $u_k$, can be cast as
\begin{equation*}
\min\limits_{u\in L^2(0,t_f;Y)} \|\mathscr{L} u-\varphi_k\|_X^2+\beta\sigma^2\int_0^{t_f}|p(t)|^2\,dt\quad 0\le k\le m
\end{equation*}
where
\begin{equation*}
p(t)=\int_t^{t_f}S^*_{s-t}C^* u(s)\, ds.
\end{equation*}
Thus, we select the parameter $\beta$ so that
$\varepsilon^2 +\beta\sigma^2$ is balanced, where $\varepsilon$ is the accuracy of the fidelity term
\begin{equation*}\mathscr{L} u_k -\varphi_k=\varepsilon.
\end{equation*}
%This provides a basic justification for the balance principle, presented in \cite{newchoice}, for regularization parameter selection.
\end{remark}

The following theorem provides the error estimate of our reconstruction method in the real 
Hilbert space setting, as well as justification for the
method based on mixed regularization.  
In short, there are two sources of error in approximating the Fourier coefficients.  The first source of error is due
to the ill-posedness of $\mathscr{L} u_k=\varphi_k$, while the second source of error is due to the noise, $\delta$, in the observed
data.  The errors must be balanced to obtain the best possible solution.  The proof is omitted, as it is a straightforward application of
the Cauchy-Schwarz inequality.

\begin{theorem}[Error Estimate]
\label{errortheorem}
Suppose $(A^*,C^*)$ is approximately controllable,
there exists $u_k\in L^2(0,t_f;Y)$ such that
\begin{equation*}
\|\mathscr{L}u_k-\varphi_k\|_X\le \varepsilon_k
\end{equation*}
for each $0\le k\le m$.
If we define,
\begin{equation*}
v(t)=\yd(t)-y(t)
\end{equation*} and
\begin{equation*}
\|v(t)\|\le\delta,
\end{equation*}
then
\begin{equation*}
\|x_0-x_\delta^m\|_X\le \|x_0-x^m\|_X + \sum\limits_{k=0}^m\left(\varepsilon_k\| x_0\|_X
+  c(\delta, t_f)\|u_k(t)\|_{Z}\right)
\end{equation*}
where $Z=L^2(0,t_f;Y)$
and
\begin{equation*}
\|x_0-x^m\|_X
\end{equation*} 
is the truncation error of the generalized Fourier series.
Furthermore, if $x_0\in C^k(\Omega)$, then
\beq
\| x_0-x^m_\delta\|_X \le \sum\limits_{k=0}^m\left(\varepsilon_k\| x_0\|_X
+ c(\delta, t_f) \|u_k(t)\|_{Z}\right)
\label{tfourier}
\eeq
for $m$ sufficiently large.
\end{theorem}

%\begin{proof}
%%\textcolor{red}{needs to be finished}
%
%%By the orhonormality of $\{\varphi_k\}_{k=0}^m$ and the Cauchy-Schwarz inequality
%%\begin{equation*}
%%\|x_0-x_\delta^m\|_X  \le \|x_0-x^m\|_X+\|x^m-x_\delta^m\|_X
%% \le \|x_0-x^m\|_X +\|\sum\limits_{k=0}^m\left(\<x_0,\varphi_k\>_X-\<u_k(t),\yd(t)\>_Z\right)\varphi_k\|_X 
%%\le \|x_0-x^m\|_X + \|\sum\limits_{k=0}^m\left(\<\mathscr{L} u_k-\varphi_k,x_0\>_X+\<u_k(t),v(t)\>_Z\right)\varphi_k\|_X
%%\end{equation*}
%%from which the estimate follows.
%%The estimate \eqref{tfourier} follows from the standard Fourier series analysis.
%Follows from standard Cauchy-Schwarz argument.
%\end{proof}
%  The error estimate is desired to be independent of the initial condition, $x_0$, however
%this is not realistic since $x_0$ is the unknown.
% It should be noted, that the error estimate is the worst possible
%error obtained in the estimation.  
Better error estimates may be realized, however, the results of Theorem \ref{errortheorem}
also provide justification for the regularization methods.  
By the estimate,
\begin{equation*}
\|x_0-x_\delta^m\|_X\le \|x_0-x^m\|_X + \|\sum\limits_{k=0}^m\left(\<\mathscr{L} u_k-\varphi_k,x_0\>_X+\<u_k(t),v(t)\>_Z\right)\varphi_k\|_X
\end{equation*}
we immediately see the need for appropriately solving $u_k$.  If the noise level, $\delta$, is large we must obtain 
controls which are sufficiently regular, so that the term
\begin{equation*}
\< u_k(t),v(t)\>_Z
\end{equation*}
is small, while simultaneously ensuring $\|\mathscr{L} u_k-\varphi_k\|_X$ is small.  The following remark further justifies the previous statement.
%\textcolor{red}{recall high frequency justification}
\begin{remark}\label{noiseremark}
Suppose the noise in the data is highly oscillatory, such as
$\cos(l\pi t)$.  Then
the error in the Fourier coefficients has the term
\beq
\int\limits_0^1 u_k(t)\cos(l\pi t)\, dt = \frac{1}{l\pi}\int\limits_0^1 u_k'(t)\sin(l\pi t)\, dt.
\label{}
\eeq
That is, the highly oscillatory parts may be damped by $l\pi$, if $u_k$ is sufficiently smooth.  
Thus, we utilize a penalty which enforces smoothness  on the control $u_k.$ 
%provides further justification for the use of a penalty term in the cost functional which enforces smoothness on the control $u_k$.
\end{remark}
  It is also apparent that
the accuracy, $\varepsilon_k$, in solving 
\begin{equation*}
\mathscr{L} u_k=\varphi_k
\end{equation*}
is necessary for an accurate reconstruction of $x_0$.  In practice, we must balance 
the accuracy of solving $\mathscr{L} u_k=\varphi_k$
and the regularity imposed on $u_k$ via the regularization methods.  This concern is addressed in Section
\ref{balancesec} where we discuss how to balance the method to obtain stable but accurate solutions.

\subsection{Variation of the Dual Control Method}
\label{alternatedual}
In this section, we outline an alternate procedure for obtaining reconstructions of the initial condition, $x_0$.  
This approach is based on the adjoint control approach developed in the previous section.  Rather than selecting a collection
$\{p_k(0)\}_{k=0}^m$ to be a basis for $X$, we
select $\{u_k(t)\}_{k=0}^m$ to be a basis(not necessarily orthonormal) for $Z=L^2(0,t_f;Y)$.  Assuming the relation 
\eqref{control} holds, we construct the adjoint set
$\{\tilde{p}_k\}_{k=0}^m$ by the relations
\begin{equation*}
\mathscr{L} u_k=\tilde{p}_k.
\end{equation*}
%where one should note that $\{\tilde{p}_k\}$ is not necessarily linearly independent.
Note that the collection $\{\tilde{p}_k\}_{k=0}^m$ is linearly independent under the assumption
that $(A,C)$ is controllable, i.e.,
\begin{equation*}
R(\mathscr{L})=X\Rightarrow N(\mathscr{L})=\emptyset.
\end{equation*}
Thus, if $(A,C)$ is exactly controllable, we form an orthogonal(orthonormal) basis
by the Gram-Schmidt method.
The coefficients of $x_0$ are computed by defining the Gram matrix
\begin{equation*}
G_{k,l}=\<\tilde{p}_k,\tilde{p}_l\>_X
\end{equation*}
and setting $\boldsymbol\beta=(\beta_0,\ldots,\beta_m)^t$ such that
\begin{equation*}
%\boldsymbol{\beta}=G^{-1}\left(\begin{array}{c}
%\int\limits_0^{t_f}\<u_0,\yd\>\, dt\\
%\vdots\\
%\int\limits_0^{t_f}\<u_m,\yd\>\, dt.
%\end{array}\right)\quad\mbox{or}\quad
\boldsymbol{\beta}=G^{-1}\left(\begin{array}{c}
\int\limits_0^{t_f}\<u_0,\xi-\yd\>\, dt\\
\vdots\\
\int\limits_0^{t_f}\<u_m,\xi-\yd\>\, dt.
\end{array}\right)
\end{equation*}
The coefficients $\boldsymbol\beta$ can be computed efficiently by the Cholesky decomposition $G=LL^*,$
since $G$ is symmetric positive definite.
Again, the algorithm is well-posed under the exact controllability \eqref{exactcontrol} of the
adjoint system which, in general,
may not be true.
If the adjoint system is not exactly controllable, care must be exercised to ensure the
set $\{\mathscr{L} u_k\}_{k=0}^m$ is linearly independent.

\begin{mdframed}[leftmargin=10pt,rightmargin=10pt]
\textbf{Variation of Dual Control Algorithm:}
\label{dualalg2}
\begin{enumerate}
\item Pick a basis $\{u_k(t)\}_{k=0}^m$ for $U_m\subset L^2(0,t_f;Y)$
\item Compute $\tilde{p}_k$ by $\mathscr{L}u_k=\tilde{p}_k$
\item Compute the Gram matrix $G_{k,l}=\<\tilde{p}_k,\tilde{p}_l\>_X$% and form the Cholesky factorization $G=LL^*$
\item Set $\boldsymbol{y}_k=\int\limits_0^{t_f}\<u_k(t),\xi(t)-y(t)\>_Y\, dt$ and compute the approximate Fourier coefficients $\boldsymbol{\beta}=G^{-1}\boldsymbol{y}$
\item Compute the approximation 
\begin{equation*}
x_0^m=\sum\limits_{k=0}^m\beta_k \mathscr{L} u_k
\end{equation*}
\end{enumerate}
\end{mdframed}
There are several potential advantages to this approach.  Namely, one can directly regulate the properties of the controls
$u_k$, such as smoothness or sparsity.  Secondly, the operator $\mathscr{L}$ does not need to be inverted.  
However, since the pair $(A^*,C^*)$ is not necessarily controllable, we are not guaranteed linear independence
of the set $\{\tilde{p}_k\}_{k=0}^m$.  Thus, solving 
\beq
G\boldsymbol\beta=\left(\begin{array}{c}
\int_0^{t_f}\<u_0,\yd\>_Y\, dt\\
\vdots\\
\int_0^{t_f}\<u_m,\yd\>_Y\, dt
\end{array}\right)
\label{}
\eeq
for $\boldsymbol\beta$ requires regularization.  This method only requires the solution of one ill-posed problem, but requires
the formation of the $m+1$ adjoints $p_k$.  Therefore, this method may be less expensive than the dual control method, however, with a tradeoff in accuracy.

\subsection{Implementation Issues}
\label{implement}
In this section, we discuss the necessary numerical issues for the implementation of the methods developed
in this section.  
%We describe the numerics for the corresponding forward problem.
For the numerical implementation for solving the dual control problem we use the 
  Crank-Nicholson scheme 
\beq
-\frac{p^{k+1}-p^k}{\Delta t}=A^*\frac{p^{k+1}+p^k}{2} + C^*u_{k+\frac{1}{2}}
\label{cn}
\eeq
for \eqref{d2}
where $u_{k+1/2}$ is evaluated at the mid-point of the interval $[t_k,t_{k+1}]$
and $t_k=kt_f\Delta t.$
At the time step $k+1$ the solution is computed by
%\beq
%p^{k+1}=-r_{1,1}(A^*\, \Delta t)p^k - \Delta t\left(I-\frac{\Delta t}{2}A^*\right)^{-1}C^*u_{k+\frac{1}{2}}
%\label{}
%\eeq
\beq
p^{k+1}=\left(I+\frac{\Delta t}{2}A^*\right)^{-1}\left(I-\frac{\Delta t}{2}A^*\right)p^k - \Delta t\left(I+\frac{\Delta t}{2}A^*\right)^{-1}C^*u_{k+\frac{1}{2}}
\label{}
\eeq
%where
%\beq
%r_{1,1}(A^*\, \Delta t)=(2-A^*\, \Delta t)^{-1}(2+A^*\, \Delta t)=\left(I-\frac{\Delta t}{2}A^*\right)^{-1}\left(I+\frac{\Delta t}{2}A^*\right).
%\label{}
%\eeq
%is the
utilizing the 
$(1,1)$ Pad\'{e} approximation for $\exp(-A^*\Delta t)$.
%The discretized RLLS inverse problem is then formulated as
%\beq
%\min\limits_{x\in X_m} \| y- M^n x\|_Y + \beta\|x\|_X
%\label{}
%\eeq
%where we have defined the controllability(observability) matrix
%\beq
%M^n=[C,Cr_{1,1}(A\, \Delta t),\ldots ,Cr_{1,1}(A\, \Delta t)^{n-1}].
%\label{}
%\eeq
%and $X_m$ is a finite-dimensional subspace of $X$.
In the dual control formulation, the discretized problem for each control $u$ is formulated 
as
\begin{equation*}
\min\limits_{u\in U_m}\|L^n u-\varphi\|_Y+\beta\psi(u)
\end{equation*}
where $L^n=(M^n)^*$, given
\beq
M^n=[C,Cr_{1,1}(A\, \Delta t),\ldots ,Cr_{1,1}(A\, \Delta t)^{n-1}]
\label{}
\eeq
and $U_m$ is a finite-dimensional subspace of $L^2(0,t_f;Y)$
and $\psi$ is a chosen penalty term.

If necessary, higher order Pad\'{e} approximations may be considered, which are of the form 
\beq
r_{m,n}(z)=\frac{P_m}{Q_n}(z)=\frac{a_0 +a_1 z+\ldots + a_m z^m}{b_0 + b_1 z+\ldots + b_n z^n}
\label{pade}
\eeq 
where the degree of $P,Q$ is not more than $m,n$ respsectively.
Higher order Pad\'{e} approximations of semigroups are discussed in detail in the paper
 \cite{stablesemigroup}.

\subsubsection*{\underline{Operator Splitting for Convection-Diffusion Equation}}
$ $\newline

\noindent
The Crank-Nicholson scheme works well for the diffusion dominant case, however, for the convection dominant case it is
necessary to solve the problem more accurately ( and such that the physics are obeyed).
In this section, we describe the numerics for the initial condition estimation of the convection-diffusion equation
\begin{equation}
\frac{\d v}{\d t}(x,t)=c(x)\cdot\nabla v(x,t) + \nabla\cdot(d(x) \nabla v)(x,t)+ f(x) ;\quad
 v(x,0)=v_0(x)
\end{equation}
where $c(x),d(x)$ are the convection and diffusion coefficients, respectively.

The reconstruction methods have a natural extension to such problems, using a differential operator splitting
\begin{equation*}
\frac{\d v}{\d t}=Lv(t)=(A+B)v(t)
\end{equation*}
where $A,B\in\mathcal{L}(X)$

For the numerical solution of the convection-diffusion equation, we consider two stage Strang operator splitting
\begin{equation*}
v(x,t+\Delta t)=S^h_{\frac{\Delta t}{2}}S^p_{\Delta t}S^h_{\frac{\Delta t}{2}}v(x,t)
\end{equation*}
where $S^p_t,S^h_t$ are the semigroups corresponding to the parabolic and hyperbolic subproblems, respectively.
That is, $S^p_t,S^h_t$ are the \cz semigroups generated by $A,B$ respectively.
 
Assuming a constant convection coefficient $c$, we solve the hyperbolic subproblem via the 
method of characteristics
$v(t_{n+1},x)=v(t_n,x-c\Delta t)$ 
where the right hand side is evaluated via cubic interpolation.
As in the previous section, we use the Crank-Nicholson method for solving the parabolic subproblem, using the approximating polynomial
\begin{equation*}
r_{1,1}(z)=\frac{2+z}{2-z}
\end{equation*}
for the approximation of $\exp(A\Delta t)$.

\section{Generalized Multi-parameter Regularization for Control Solution}
This section is devoted to discussing the multi-term regularization method utilized for solving \eqref{control}, without going into detail.
In general, rather than solving
$$
\mathscr{L}u=\varphi
$$
directly, we seek a minimum of
\beq
\mathcal{J}_\beta(u) = \phi(u,\varphi)+\beta\psi(u),
\label{gentik}
\eeq
over $u\in C$,
where the fidelity term, $\phi$, is chosen based on the noise statistic, while $\psi$ is chosen based on which class the solution $x$ should belong to.
Whenever $\phi(u,\varphi)=\| \mathscr{L}u-\varphi\|_X^2$ and $\psi(u)=\| u\|_X^2$ this formulation coincides with the classical Tikhonov regularization
$$
\min\limits_{u\in\mathcal{C}}\frac{1}{2}\|\mathscr{L}u-\varphi\|_X^2+\frac{\beta}{2}\|u\|_X^2.
$$
The main drawback to this method is the single regularization term $\psi.$
Modern day scientific problems typically involve applications where the standard Tikhonov regularization fails to capture the full set of distinct
features in the physical solution.  
Many research efforts have been devoted to improving the standard regularization techniques for a wide range of applications (see \cite{lions,itokuna,perona,rof} for example).
It is not a goal of this paper to cover this in detail.  These references and the references therein provide a thorough study of such methods.
Especially in the field of image processing, the solution often exhibits a multiscale structure typically described by multi-resolution analysis.
In such applications, single parameter regularization can oversmooth the solution such as the case of $\psi=\|\cdot\|_{L^2}^2$ or exhibit stair-case effects such as the case of
$\psi=\|\cdot\|_{TV}.$  In order to capture the multiscale structure of solutions without introducing oversmoothing or staircasing, many research efforts have
focused on mixed regularization approaches, such as combining the $L^2$ penalty term with the $TV$ penalty:
\beq
\min\limits_{u\in C} \frac{1}{2}\int\limits_{\Omega}|\mathscr{L}u-\varphi|^2\, d\xi+\frac{\eta_1}{2}\int\limits_\Omega |u|^2\, d\xi+\eta_2\int\limits_\Omega |\nabla u|\, d\xi.\label{mixedT}
\eeq

To capture multi-scale solution profiles, we employ the multi-parameter Tikhonov regularization technique, i.e., we minimize
\begin{equation}
\mathcal{J}_{\boldsymbol\eta}(u)=\phi(u,\varphi)+\boldsymbol\eta\cdot\boldsymbol\psi(u).
\label{multitik}
\end{equation}
The terms $\phi,\boldsymbol\psi$ are known as the fidelity and regularization terms, respectively.  Here, $\{\psi_k\}_{k=1}^n$ is the set of
regularization terms, $\{\eta_k\}_{k=1}^n$ are the regularization parameters, and we take the dot product
 \begin{equation*}
 \boldsymbol\eta\cdot\boldsymbol\psi(u)=\sum\limits_{k=1}^n\eta_k\psi_k(u)
 \end{equation*}
for $\boldsymbol\eta=(\eta_1,\eta_2,\ldots,\eta_n)$ and $\boldsymbol\psi(u)=(\psi_1(u),\psi_2(u),\ldots,\psi_n(u)).$
The functionals $\phi,\boldsymbol\psi$ can be chosen based on any a priori information about the problem and its exact solution.
Then,
\begin{equation*}
u_{\boldsymbol\eta}=\argmin{u\in C} \,\,\mathcal{J}_{\boldsymbol\eta}(u)
\end{equation*}
 is taken as the regularized solution.
For instance, in the case of a multiscale image with a smooth region and a stepped region, one may consider the $L^2$-$TV$ regularization \eqref{mixedT}.
In this work, we also consider the penalty term 
$$
\psi(u)=\|u\|_{\ell_p}^p\quad p\le 1,
$$
to enforce sparsity in the solution.

%A common theme in the literature and throughout this paper, is the introduction of sparsity.  
%Due to the ever increasing size of today's most challenging and interesting problems, it is often necessary 
%to reduce the size of a problem and/or obtain a solution with as few nonzeros as possible.  
%This idea corresponds to minimizing the size of a dataset or solution by keeping only what information 
%is necessary in order to describe the solution accurately.  In the case of obtaining solutions with few nonzeros, we consider regularization terms such as
%\begin{equation*}
%\psi_k(x)=\|x\|^p_{\ell_p}:=\sum\limits_{k=1}^\infty |x_k|^p\quad\mbox{for } p\le 1.
%\end{equation*}
%This particular choice of $\psi$ requires the treatment of the associated numerical difficulties.

\subsection{Balance principle}
\label{balancesec}
We discuss here the balance principle for the single-term regularization, based on the paper \cite{newchoice, augtik}.
Prior to this selection rule, most selection rules (e.g. Morozov's discrepancy principle) were based on either the performance level (noise) 
$$
\phi(u,\varphi)
$$
or the complexity level
$$
\psi(u,\varphi),
$$
alone.  The selection rule
developed in \cite{newchoice, augtik} is based on balancing the performance level and the complexity level.
Consider maximizing the conditional density $p((u,\tau,\lambda)|\varphi)\sim p(\varphi|(u,\tau,\lambda))p(u,\tau,\lambda)$ where $(\tau,\lambda)$ are 
density functions for $\phi,\psi$, respectively, both having Gamma distribution.  The balancing principle is derived
from the Bayesian inference \cite{augtik}
\begin{equation}
\min\limits_{(u,\tau,\lambda)} \tau\phi(u,\varphi)+\lambda\psi(u)+\tilde{\eta}_0\lambda-\tilde{\alpha_0}\ln\lambda +\tilde{\eta}_1\tau-\tilde{\alpha}_1\ln\tau.
\label{bayesfunct}
\end{equation}
By definition, $(u,\lambda^*,\tau^*)\in X\times\mathbb{R}^+\times\mathbb{R}^+$ is a critical point of \eqref{bayesfunct} if
\begin{equation*}
\begin{array}{l}
u^* =\argmin{u}\{\phi(u,\varphi)+\lambda^*(\tau^*){-1}\psi(u)\}\\
\\
\psi(u^*)+\eta_0-\alpha_0\frac{1}{\lambda^*}=0\\
\\
\phi(u^*,\varphi)+\eta_1-\alpha_1\frac{1}{\tau^*}=0
\end{array}
\end{equation*}
By optimality, the regularization parameter satisfies
\begin{equation}
%u_\eta  = \argmin{u}\{\phi(u,\varphi)+\eta\psi(u)\},\quad
\eta^*  = \frac{\alpha_0}{\alpha_1}\frac{\phi(u^*,\varphi)+\eta_1}{\psi(u^*)+\eta_0}.
%,\quad \mu=\frac{\tilde{\alpha}_1}{\tilde{\alpha}_0.}
\label{balanceprinc}
\end{equation}
The authors arrive at the selection criterion
$$
\eta^*=\frac{\alpha_0}{\alpha_1^{1-d}}\frac{\phi(u^*,\varphi)^{1-d}}{\psi(u^*)+\eta_0},\quad 0< d<1,
$$
by rescaling $\alpha_0$ as $\sigma_0^{-d}$, in order to ensure the conditions
$$
\lim\limits_{\sigma_0\to 0}\bar{\eta}(\sigma_0)=0,\quad\lim\limits_{\sigma_0\to 0}\frac{\sigma_0^2}{\bar\eta(\sigma_0)}=0
$$
are satisfied, where $\sigma_0^2$ is the variance.
%The Bayesian inference selection corresponds to minimization of the functional
%\beq
%\frac{F(u_\eta)^{1+c}}{\eta.}
%\label{}
%\eeq
%The minimum occurs when $\frac{d}{d\eta}F(u_\eta)=0$ so that
%\beq
%\frac{F(u_\eta)^{c}\left((1+c)F'(u_\eta)\eta-F(u_\eta)\right)}{\eta^2}=0
%\label{}
%\eeq
%which implies
%\beq
%\phi(u_\eta,\varphi)=c\eta\psi(u_\eta)
%\label{balance}
%\eeq
%since $F'(u_\eta)=\psi(u_\eta).$
%The relation (\ref{balance}) is what we refer to as the balance principle for choosing the parameter $\eta$.
Further discussion on the validity of this method, as well as the selection of the constants $\alpha_i,\eta_0,d$, can be found
in \cite{newchoice}.
The following iterative algorithm for determining $u_\eta,\eta$ is utilized: 
\begin{mdframed}[leftmargin=10pt,rightmargin=10pt]
\textbf{Iterative algorithm to solve for $(u_\eta,\eta$):}\label{choicerule}

\noindent
Choose an initial guess $\beta_0 > 0$, and set $k = 0$. Find $(u_k, \beta_k)$ for $k\ge 1$ as follows:
\begin{enumerate}
\item \label{step1} Solve for $u_{k+1}$ by the Tikhonov regularization method to obtain
$$
u_{k+1}=\argmin{u}=\{\phi(u,\varphi)+\beta_k\psi(u)\}.
$$
\item Update the regularization parameter $\beta_{k+1}$ by
$$
\beta_{k+1}=\alpha\frac{\phi(u_{k+1},\varphi)^{1-d}}{\psi(u_{k+1})+\eta_0}.
$$
\item If a stopping criterion is met, stop, else set $k=k+1$ and repeat from step \ref{step1}.
\end{enumerate}
\end{mdframed}

%The natural choice for updating the parameter $\eta$ is the fixed point iterate
%\begin{equation}
%\eta^+  =  \frac{1}{\gamma}\frac{\phi(u_\eta,\varphi)}{\psi(u_\eta)} \label{}
%\end{equation}
%where $u_\eta$ is the solution with the previous value for $\eta$.
%Here, $\gamma$ is selected by the two-step procedure proposed in \cite{regparam,newchoice} :
%\begin{itemize}
%\item Choose $\gamma_0,\eta_0$
%\item Compute $u_\eta$ with chosen parameters
%\item Set $\gamma$ by
%\beq
%\gamma=\gamma_0\left(\frac{\phi(u_\eta,\varphi)}{.05\phi(0,\varphi)}\right)^d
%\label{}
%\eeq
%\end{itemize}
%with $d,\gamma_0$ heuristically chosen as $d=\frac{1}{4},\gamma_0=10$.

%The multi-parameter balance principle is obtained by minimizing
%\begin{equation*}
%\Phi_\gamma(\boldsymbol\eta)=\frac{\gamma^\gamma}{(\gamma+2)^{\gamma+2}}\frac{F^{2+\gamma}(\boldsymbol\eta)}{\underset{i}{\Pi}\eta_i}
%\end{equation*}
%where
%\begin{equation*}
%F(\boldsymbol\eta)=\mathcal{J}_{\boldsymbol\eta}(x_{\boldsymbol\eta})
%\end{equation*}
%with $x_{\boldsymbol\eta}$ the minimizer of \eqref{multitik}.

\section{Numerical Tests}
\label{applications}

\subsection{1-D Diffusion Equation}
In this section, we consider inverse problems involving the 1-D diffusion equation
\begin{equation}
\begin{array}{c c}
& \frac{\d v}{\d t}=\frac{\d}{\d x}(d(x)\frac{\d v}{\d x}) +f(x)\quad x\in \Omega\subset \mathbb{R} \\
& \\
& v(0,t)=0=v(1,t)\\
& \\
& y(t)=Cv(t)
\label{1ddiffex}
\end{array}
\end{equation}
with Dirichlet boundary conditions,
where the measurements are restricted to a subinterval $\Omega_s\subset\Omega=[0,1]$, for the time $0\le t\le 1$.
%The measurements are taken as averages over two patches $\Omega_1,\Omega_2$.
Specifically, the operator $C$ takes average measurements over the two intervals $(.23,.31)$ and $(.46,.53).$
The thermal conductivity, $d$, is potentially variable in space, but known.
The 1-D diffusion equation is formulated as an abstract Cauchy problem \eqref{cauchy}
where
$$
Av=\frac{d }{dx}(d(x)\frac{dv}{dx})
$$
and
\begin{equation*}
\begin{array}{c c}
\mbox{dom}(A) =\{& v\in L^2(0,1)| v,\frac{dv}{dx}\mbox{ are absolutely continuous,}\\
& \frac{d^2 v}{d x^2}\in L^2(0,1)\mbox{ and } v(0)=0=v(1)\}.
\end{array}
\end{equation*}
It is a standard exercise to show that $A$ generates a strongly continuous semigroup.

\subsection*{Example 1 : Spatially varying diffusion coefficient}
\begin{figure}[!hh]
\centering
\subfloat[][]{\label{1dvari-a}\includegraphics[width=3in]{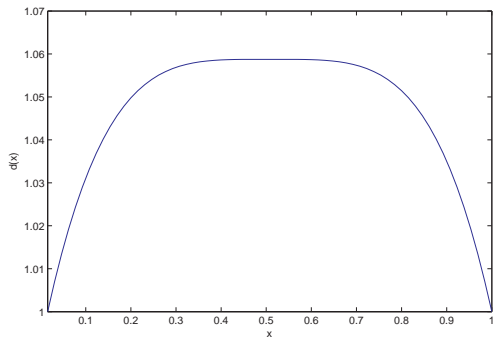}}\quad
\subfloat[][]{\label{1dvari-b}\includegraphics[width=3in]{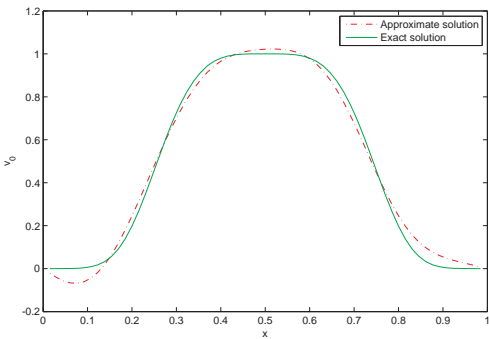} }
\caption[1-D reconstruction for variable coefficient diffusion equation]{\subref{1dvari-a} Thermal conductivity; 
\subref{1dvari-b} Reconstruction with 10\% noise in measurements versus Exact initial condition.}
 \label{1dvari}
\end{figure}

For this example, we consider the case when the thermal conductivity is spatially variable.  In particular, we
take
$$
d(x) = 1.0625-(x-\frac{1}{2})^4
$$
and the initial condition is given by
$$
v_0(x)=e^{-200(x-\frac{1}{2})^4}.
$$
%No filter ...compare coiflet/versus standard basis using same diffusion...multiply d by .01
We take the basis $\varphi_k=\{\sin(k\pi x)\}_{k=0}^m$ and
solve for the controls, $u_k,$ using the $L^1$-$H^1$ regularization method, for $m=8$.
It should be pointed out that the abstract Cauchy based dual control method does not make any assumptions
on the coefficients of the PDE.  
 Assuming a noise level of 10\%, we obtain an accurate and stable reconstruction
with the parameters $\eta_1=5\times 10^{-8},\eta_2=1\times 10^{-10}$.  
The corresponding results are depicted in Figures \ref{1dvari}.

\subsubsection*{\underline{Comparison of basis choices}}
\begin{figure}[!hh]
\centering
%\subfigure[][]{\label{1dvari2-a}
\includegraphics[width=3in]{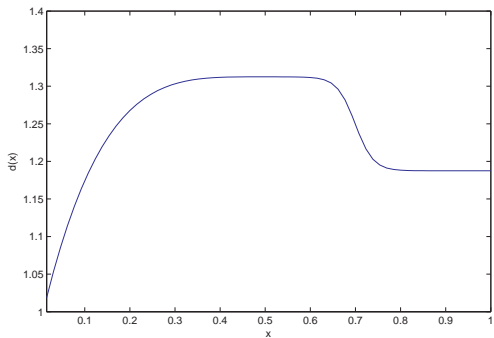}
%\subfigure[][]{\label{1dvari2-b}\includegraphics[width=3in]{/home/crhumber/ncsuthesis-0.2.2/results/heat1d/1dvari3.eps} }
\caption[1-D reconstruction for variable coefficient diffusion equation]{ Thermal conductivity.}
 \label{1dvari2}
\end{figure}

\begin{figure}[!hh]
\centering
%\subfloat[][]{\label{1dvari3-a}\includegraphics[width=3in]{/home/crhumber/ncsuthesis-0.2.2/results/heat1d/daub5e71e10.eps}}\quad
%\subfloat[][]{\label{1dvari3-b}\includegraphics[width=3in]{/home/crhumber/ncsuthesis-0.2.2/results/heat1d/sine5e71e10.eps} }
\subfloat[][]{\label{1dvari3-a}\includegraphics[width=3in]{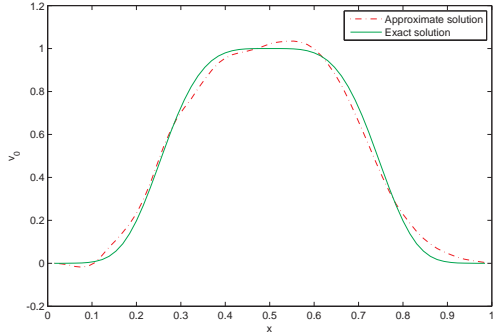}}\quad
\subfloat[][]{\label{1dvari3-b}\includegraphics[width=3in]{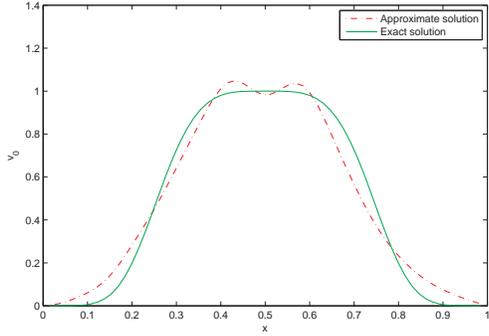} }
\caption[Comparison using different basis functions]{\subref{1dvari3-a} Reconstructed initial condition using Daubechies-10 wavelets
with $\eta_1=5\times 10^{-7},\eta_2= 1\times 10^{-11}$; 
\subref{1dvari3-b} Reconstructed initial condition using sine basis with $\eta_1=5\times 10^{-6},\eta_2= 1\times 10^{-10}$.}
 \label{1dvari3}
\end{figure}

Here, we compare the reconstructions obtained by two different basis choices.
For this example, we
take
$$
d(x) = \left\{\begin{array}{l l}
1\frac{5}{16}-5(x-\frac{1}{2})^4,& \quad 0\le x< \frac{1}{2}\\
\\
1\frac{3}{16}+\frac{1}{8+e^{-50(x-.65)}},& \quad\frac{1}{2}\le x\le 1,
\end{array}\right.
$$
as depicted in Figure \ref{1dvari2},
and the initial condition is given by
$$
v_0(x)=e^{-200(x-\frac{1}{2})^4}.
$$
 We
solve for the controls $u_k$ using the $L^1$-$H^1$ regularization method.
 Assuming a relative noise level of 5\%, we obtain an accurate and stable reconstruction
with the parameters $\eta_1=5\times 10^{-7},\eta_2=1\times 10^{-11},$ by computing $m=8$ coefficients,
using Daubechies wavelets.  
%The corresponding results are depicted in Figures \ref{1dvari2}. 
In Figure \ref{1dvari3}, one can see a comparison of two reconstructions using a standard sine basis and Daubechies-10
 wavelets.  The reconstruction obtained
using the Daubechies wavelets is more accurate and stable than the sine basis reconstruction, even with well-tuned 
regularization parameters.  This example illustrates how the basis choice affects the resulting
reconstruction.

%\clearpage
\subsection{2-D Diffusion Equation}
In this section, we consider inverse problems involving the 2-D diffusion equation
\begin{equation} \label{2ddiff}
\begin{array}{c c}
& \frac{\d v}{\d t}(\boldsymbol{x},t)=d\Delta v(\boldsymbol{x},t) + f(\boldsymbol{x})\\
& \nonumber \\
& v(\boldsymbol{x},t)=0\quad x\in\d\Omega\\
& \nonumber \\
& v(\boldsymbol{x},0)=v_0(\boldsymbol{x})
\end{array}
\end{equation}
where
$\boldsymbol{x}=(x,y)\in\Omega\subset\mathbb{R}^2.$
As in the 1-D case, we work on the time interval $0\le t\le 1.$

The 2-D diffusion equation can be cast in the abstract Cauchy framework
where
$A$ coincides with the closure of the Laplace operator, defined by
$$
\Delta f(\boldsymbol{x}) := \frac{\d^2 f}{\d x^2}+\frac{\d^2 f}{\d y^2} 
$$
for every $f$ in the Schwartz space

$$
\mathscr{S}(\mathbb{R}^n):=\left\{f\in C^\infty(\mathbb{R}^n): \lim\limits_{|x|\to\infty}|x|^kD^\alpha f(x) =0 \mbox{ for all } k\in \mathbb{N}\mbox{ and }\alpha\in\mathbb{N}^n\right\}.
$$

\subsection{2-D Convection-Diffusion Equation}
In this section, we present severeal numerical results for inverse problems involving the convection-diffusion equation
\begin{equation}
\begin{array}{c c}
& \frac{\d v}{\d t}=c(x)\cdot\nabla v + \nabla\cdot(d(x) \nabla v)+ f(x) \\
& v(x,0)=v_0(x).
\end{array}
\end{equation}
For the results presented here, we assume $c(x)\equiv c, d(x)\equiv d$ are constant (or at least locally constant), and we take $f(x)\equiv 0$.  For both simulations, the domain is taken as the unit square $\Omega=[0,1]\times[0,1].$

\subsection*{Initial condition reconstruction}
\begin{figure}[!hh]
\centering
\subfloat[][]{\label{2dcdmeas-a}\includegraphics[width=3in]{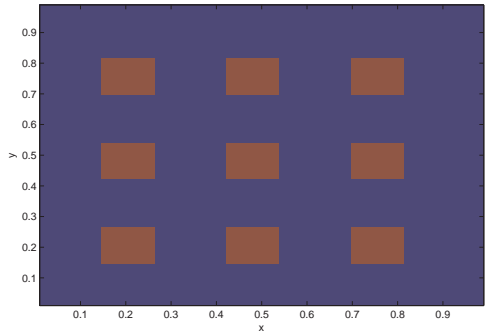}}\quad
\subfloat[][]{\label{2dcdmeas-b}\includegraphics[width=3in]{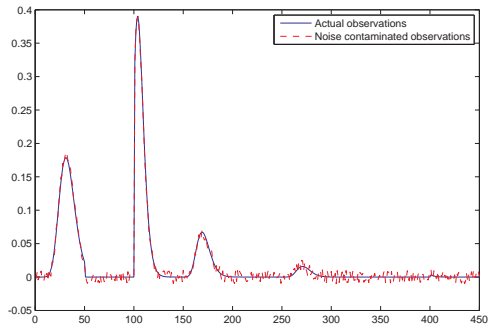} }
\caption[Measurements and locations for 2-D convection-diffusion equation]{\subref{2dcdmeas-a} Nine measurement locations depicted in red; 
\subref{2dcdmeas-b} Noisy measurements used for reconstruction compared with exact measurements.}
 \label{2dcdmeas}
\end{figure}

\begin{figure}[!hh]
\centering
\subfloat[][]{\label{2dcd-a}\includegraphics[width=3in]{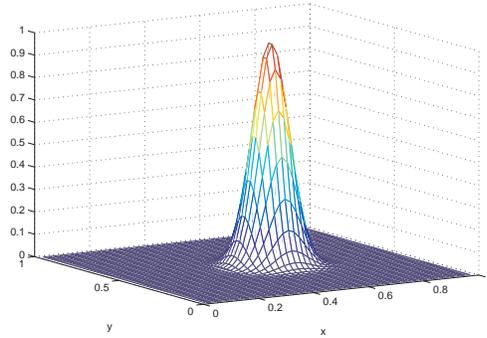}}\quad
\subfloat[][]{\label{2dcd-b}\includegraphics[width=3in]{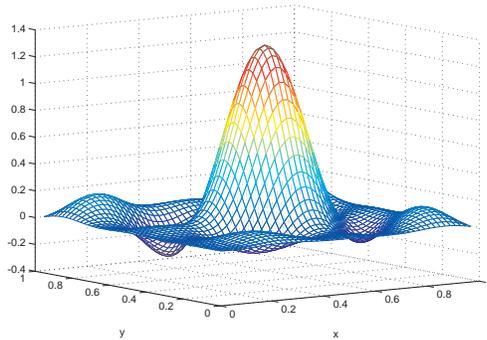} }
\caption[Reconstruction for 2-D convection-diffusion equation]{\subref{2dcd-a} Exact initial condition; 
\subref{2dcd-b} Reconstruction with 10\% noise in measurements.}
 \label{2dcd}
\end{figure}

We consider the initial condition reconstruction problem with $d=.1,c=(\frac{1}{2},\frac{1}{2})$ known, where the observation operator is defined by
$$
Cv(t)=\frac{1}{\mu(\Omega_s)}\int\limits_{\Omega_s}v(s) d\mu,  
$$
where $\mu(\Omega_s)$ is the volume of the set $\Omega_s.$
That is, we take average measurements over a sample set $\Omega_s\subset \Omega$.  For this simulation, we take nine
measurement locations equally spaced over the domain, each location of size $\frac{1}{10}\times\frac{1}{10}$.  
The corresponding contaminated measurements are depicted in Figure \ref{2dcdmeas-b}. 
%which can filtered via an $H^1$ (smoothing) filter prior to implementing the method.
Using the operator splitting technique outline in Section \ref{implement}, the convection-diffusion equation fits into the abstract framework \eqref{cauchy}.
The exact initial condition is
$$
v_0(x,y)=e^{-100((x-.55)^2+(y-.5)^2)}
$$
and we solve the corresponding inverse problem using the $L_1$-$H_1$ regularization, with basis functions
$$
\varphi_{k,l}(x,y)=\sin(k\pi x)\sin(l\pi y).
$$
As can be seen by comparing Figures \ref{2dcd-a} and \ref{2dcd-b}, the method for reconstructing the initial condition performs well
with the parameters $\eta_1=.03,\eta_1=1\times10^{-8}$. 
Depending on the basis choice, small errors are expected due to the truncation of the generalized Fourier series. 
In this case, we have small oscillations indicative of the sinusoidal basis.
%We also consider the convection-diffusion equation to investigate how adding convection
%affects the estimation results.  We compare the reconstruction of the same initial condition from the
%convection-diffusion equation and the diffusion equation using the same measurement locations.
%As can be seen in Figure \ref{2dcompare}, the addition of convection has the potential to sharpen results depending on measurement
%location.  In this case, there is the appearance of a hump over one measurement location when no convection is present.
%This is likely due to difficulty in capturing the direction in which information(heat flow) is propagating.

\section{Concluding remarks}

The abstract Cauchy problem provides a unified framework for the analysis of systems governed
by PDE.  The methods developed in this paper allow for the systematic reconstruction of initial
conditions of the abstract Cauchy problem.  In particular, the dual control method coupled with the
multi-parameter regularization yields a method that is very tunable and robust.  
By an appropriate basis selection for the problem at hand, and by selecting the parameters in the regularization
framework based on the balance principle, a reconstruction filter is determined based on the governing PDE.
Depending on the problem size, there may be significant overhead in computing the controls \eqref{control}.  
However, once computed, the controls can be banked (or stored) for future use.
Thus, if one carefully selects the basis and the parameters are tuned to the noise and a priori
information about the solution, 
the method can potentially be implemented
in real time, simply by integrating the controls against the data.

Diffusion processes and parabolic equations fit particularly well into this framework, due to
the necessity for stabilizing the dynamics backward in time.  The method accurately reconstructs both
initial conditions and point sources of diffusion processes, and allows the forecasting of future states.
Thus, the tool provided is valuable for problems where numerous calculations are required based on sensor
data, and for problems where integrating forward and backward in time is important.

Based on the multi-parameter regularization, the methods developed are particularly suited
for problems involving a locally supported source, such as point sources, as well as those with
sparsely distributed data.  The sparsity optimization works well for both identifying initial conditions/sources 
that are locally supported, as well as for selecting the necessary control profile.

Certain questions still remain and extensions to more difficult problems can be realized.
Specifically, nonlinear problems can be treated in a similar manner, through the development
of nonlinear dual control filters.  In a forthcoming paper, we describe the nonlinear method
for equations such as the one-dimensional viscous Burger's equation
$$
u_t+\left(\frac{u^2}{2}\right)_x=\varepsilon u_{xx},
$$
the Korteweg-de Vries (KdV) equation
$$
u_t-6uu_x=u_{xxx},
$$
and its generalizations (e.g. the Novikov-Veselov equation).
We are also interested in inverse problems regarding the incompressible Navier-Stokes equations
\begin{eqnarray*}
\frac{\d}{\d t}v_i + \sum\limits_{j=1}^nv_j\frac{\d v_i}{\d x_j} &  = & \nu\Delta v_i -\frac{\d p}{\d x_i} + f_i(x,t)\quad (x\in \mathbb{R}^n,t\ge 0),\\
& & \mbox{div}\, v = \sum\limits_{i=1}^n\frac{\d v_i}{\d x_i}=0\quad (x\in \mathbb{R}^n,t\ge 0),
\end{eqnarray*}
with initial conditions
$$
v(x,0)=v_0(x)\quad (x\in\mathbb{R}^n).
$$
The three-dimensional Navier-Stokes equations have important applications, such as weather modeling, aircraft design, rheology, etc.
Due to the practical need for considering three-dimensional Navier-Stokes, efficiency must be addressed.
%The consideration of the Navier-Stokes equations also involves the necessity for considering
%higher-dimensional problems.
In this case,
the solution for the controls must be performed efficiently, though, once computed, this framework
may be ideal for large scale problems since the filter coefficients can simply be banked.
Thus, future research for this method also involves addressing computational efficiency.

\appendix

\medskip

Received September ...

% Please replace the following email addresses by yours.

\medskip
 {\it E-mail address: }cary.humber@navy.mil\\
 \indent{\it E-mail address: }kito@math.ncsu.edu\\

\end{document}